\documentclass{elsart}

\usepackage{amsmath}
\usepackage{amssymb}
\usepackage[all]{xy}
\usepackage{sty1}

\journal{Journal of Functional Analysis}

\makeatletter
\@addtoreset{equation}{section}
\makeatother

\def\cstarg{\mathrm{C}^*(G)}
\def\fal#1{\mathrm{A}(#1)}
\def\falg{\mathrm{A}(G)}
\def\falh{\mathrm{A}(H)}
\def\falgh{\mathrm{A}(G\cross H)}
\def\falbg#1{\mathrm{A}_{#1}(G)}
\def\falbh#1{\mathrm{A}_{#1}(H)}
\def\falkg{\mathrm{A}_K(G)}
\def\falskg{\mathrm{A}_{sK}(G)}

\def\fsalg{\mathrm{B}(G)}
\def\fsalgmn{\mathrm{B}(G/N)}
\def\fsalgn{\mathrm{B}(G\! :\! N)}
\def\falcg{\mathrm{A}_c(G)}
\def\falh{\mathrm{A}(H)}
\def\blone#1{\mathrm{L}^1(#1)}
\def\bloneg{\mathrm{L}^1(G)}
\def\bltwo#1{\mathrm{L}^2(#1)}
\def\bltwog{\mathrm{L}^2(G)}
\def\id{\mathrm{id}}
\def\lebfalg{\mathrm{LA}(G)}
\def\mtg{\mathrm{MT}(G)}
\def\optens{\hat{\otimes}}
\def\ptens{\otimes^\gamma}
\def\stens{\check{\otimes}}
\def\ran{\mathrm{ran}\,}
\def\seg{\mathrm{S}}
\def\sfalg{\mathrm{S}\mathrm{A}(G)}
\def\sofalg{\mathrm{S}_0(G)}
\def\sofalh{\mathrm{S}_0(H)}
\def\sofalgh{\mathrm{S}_0(G\cross H)}
\def\sofalb#1{\mathrm{S}_0(#1)}
\def\soneg{\mathrm{S}^1(G)}
\def\tg{\mathrm{T}(G)}
\def\tk{\mathrm{T}(K)}
\def\tb#1{\mathrm{T}(#1)}
\def\vng{\mathrm{VN}(G)}
\def\vnh{\mathrm{VN}(H)}

\def\wtwog{\mathrm{W}^2(G)}

\def\supp{\mathrm{supp}}

\def\mat#1#2{\mathrm{M}_{#1}(#2)}
\def\matm#1{\mathrm{M}_m(#1)}
\def\matn#1{\mathrm{M}_n(#1)}

\def\smatn{\mathrm{M}_n}
\def\cbop#1{\mathcal{CB}(#1)}
\def\stmatn{\mathrm{T}_n}
\def\tmatn#1{\mathrm{T}_n(#1)}
\def\tmatm#1{\mathrm{T}_m(#1)}
\def\tmatnm#1{\mathrm{T}_{nm}(#1)}
\def\tmatinfty#1{\mathrm{T}_\infty(#1)}
\def\tmatfin#1{\mathrm{T}_{\mathrm{fin}}(#1)}
\def\stmatinfty{\mathrm{T}_\infty}

\def\blonenorm#1{\left\|#1\right\|_{\mathrm{L}^1}}
\def\cbnorm#1{\left\|#1\right\|_{cb}}
\def\fnorm#1{\left\|#1\right\|_{\mathrm{A}}}

\def\bltwonormb#1#2{\left\|#1\right\|_{\mathrm{L}^2(#2)}}
\def\sfnorm#1{\left\|#1\right\|_{\mathrm{S}\mathrm{A}}}

\def\sonenorm#1{\left\|#1\right\|_{\mathrm{S}^1}}

\newtheorem{tmatop}{Proposition}[section]
\newtheorem{tmatop1}[tmatop]{Corollary}
\newtheorem{standard}[tmatop]{Proposition}
\newtheorem{localisation}[tmatop]{Corollary}
\newtheorem{localisation1}[tmatop]{Corollary}

\newtheorem{feichdesc}{Lemma}[section]
\newtheorem{maintheorem}[feichdesc]{Theorem}
\newtheorem{maintheorem1}[feichdesc]{Corollary}
\newtheorem{maintheorem2}[feichdesc]{Corollary}
\newtheorem{maintheorem3}[feichdesc]{Corollary}
\newtheorem{maintheorem4}[feichdesc]{Corollary}

\newtheorem{tensorprod}{Theorem}[section]
\newtheorem{tensorprod1}[tensorprod]{Corollary}
\newtheorem{restrict}[tensorprod]{Theorem}
\newtheorem{tmultipliers}[tensorprod]{Theorem}
\newtheorem{ideals}[tensorprod]{Proposition}
\newtheorem{averaging}[tensorprod]{Theorem}
\newtheorem{isomorphism}[tensorprod]{Theorem}

\begin{document}
\begin{frontmatter}

\title
{Operator space structure on Feichtinger's Segal algebra}

\author{Nico Spronk}\ead{nspronk@math.uwaterloo.ca}


\thanks{The author's research supported by NSERC under grant 
no.\ 312515-05.}

\address{Department of Pure Mathematics, University of Waterloo, \\ Waterloo 
ON, N2L 3G1, Canada}

\bigskip
\begin{center}
{\it In memory of my friend, Sean Crawford Andrew.}
\end{center}

\begin{abstract}
We extend the definition, from the class of abelian groups to a general locally
compact group $G$, of Feichtinger's remarkable Segal algebra $\sofalg$.  
In order to obtain functorial properties for non-abelian groups, 
in particular a tensor product formula,
we endow $\sofalg$ with an operator space structure.  With this structure
$\sofalg$ is simultaneously an {\it operator Segal algebra} of the Fourier 
algebra $\falg$, and of the group algebra $\bloneg$.  We show that this 
operator space structure is consistent with the major functorial properties: 
(i)
$\sofalg\optens\sofalh\cong\sofalgh$ completely isomorphically
(operator projective tensor product), if $H$ is another locally compact group; 
(ii) the restriction map $u\mapsto u|_H:\sofalg\to\sofalh$ is completely 
surjective, if $H$ is a closed subgroup;
and (iii) $\tau_N:\sofalg\to\sofalb{G/N}$  is completely surjective, where $N$
is a normal subgroup and $\tau_Nu(sN)=\int_Nu(sn)dn$.
We also show that $\sofalg$ is an invariant for $G$ when it is treated
simultaneously as a pointwise algebra and a convolutive algebra.
\end{abstract}

\maketitle

\begin{keyword}
Fourier algebra, Segal algebra, operator space.
\MSC Primary 43A30, 43A20, 46L07;
Secondary 47L25, 46M05.
\end{keyword}

\end{frontmatter}


\section{Introduction and Notation}

\subsection{History}
In \cite{feichtinger}, Feichtinger defined, for any abelian group $G$, 
a Segal algebra $\sofalg$ of $\bloneg$.  This Segal algebra is the minimal 
Segal algebra in $\bloneg$ which is closed under pointwise multiplication by 
characters and on which multiplication by any character is an isometry.
It is proved in \cite{feichtinger}, that the Fourier 
transform induces an isomorphism $\sofalg\cong\sofalb{\hat{G}}$ where $\hat{G}$
is the dual group.  Thus, we also have that $\sofalg$ is a Segal algebra in
the Fourier algebra $\falg$ $(\cong\blone{\hat{G}})$, i.e.\ a dense ideal of 
$\falg$ which has a norm under which it is a Banach
$\falg$-module.  In fact, it is the minimal Segal algebra in $\falg$ which is
translation invariant and on which translations are isometries.

For a general, not necessarily abelian, locally compact group $G$, the Fourier
algebra is defined by Eymard \cite{eymard}.  There are hints in 
\cite{feichtinger81} of how to define $\sofalg$, as a Segal algebra in $\falg$.
We develop this fully.  $\sofalg$ is also a Segal algebra in the classical 
sense, i.e.\ a Segal algebra of $\bloneg$.
We also develop, for general locally compact groups,
the functorial properties which Feichtinger proved for
abelian groups \cite[Thm.\ 7]{feichtinger}.  One of Feichtinger's results  
is a tensor
product formula:  if $G$ and $H$ are locally compact abelian groups, then
there is a natural isomorphism $\sofalg\ptens\sofalh\cong\sofalgh$ (projective 
tensor product).  For non-abelian $G$ and $H$ we cannot expect that
$\falg\ptens\falh\cong\falgh$, in general, by \cite{losert}.  The theory
of operator spaces, and the associated operator tensor product, allows
us to obtain a satisfactory result from \cite{effrosr}:  
$\falg\optens\falh\cong\falgh$.  Thus we
are motivated to find a natural operator space structure on $\sofalg$, for a
general locally compact $G$, which allows us to recover a tensor product
formula.  Analogous to the fact that $\sofalg$ does not have a fixed natural 
norm, but rather a family of equivalent norms, we find that our operator
space structure is determined only up to complete isomorphism.
In order to deal with our operator space structure, we find it more
natural to deal with certain ``dual'' type matrices $\tmatn{\fV}$ over an 
operator space $\fV$, as opposed to the usual matrices $\matn{\fV}$, which 
forces us to summarise a coherent theory of these in 
Section \ref{ssec:opspace}.

To underscore the naturality of our operator space structure we examine
the other two major functorial properties, restriction to a closed subgroup
and averaging over a closed normal subgroup.  We show that our operator space 
structure is natural in the sense that it gives complete surjections, onto
the Feichtinger algebra of the closed subgroup in the case of restriction, and
onto the Feichtinger algebra of the quotient group in the case of averaging.
See Sections \ref{ssec:restriction} and \ref{ssec:averaging}. 

In \ref{ssec:isomorphism} we discuss an isomorphism theorem, characterising
$\sofalg$ and an invariant of $G$.  This makes no use of our operator space
structure.  As we note, our result is not even specific to $\sofalg$,
but can be applied to many spaces which are simultaneously Segal algebras
in $\bloneg$ and in some regular Banach subalgebra $\fA$ of functions on
$G$, having Gelfand spectrum $G$.

The Segal algebra $\sofalg$ seems interesting in and of itself simply for
its wealth of structure and functorial properties.  However, in the abelian
case, $\sofalg$ is a fundamental example of {\it Wiener amalgam spaces}
and of {\it modulation spaces}, which appear to be of tremendous use in
time-frequency analysis.  See \cite{feichtingerg} and references therein.
We hope our $\sofalg$, for non-abelian $G$, may prove as inspirational and
useful.

\medskip
The author is grateful to Ebrahim Samei and Antoine Derighetti, for indicating 
to him the article \cite{lohoue}, and helping to clarify some of the
results therein.  Though \cite{lohoue} was not ostensibly used, in the end,
it inspired the author to find the proof of the result of Section
\ref{ssec:averaging}.

\subsection{Harmonic Analysis}  Let $G$ be a locally compact group
with fixed left Haar measure $m$.  We will denote integration of a function $f$
with respect to $m$ variously by $\int_G fdm$ or $\int_Gf(s)ds$. 
For any $1\leq p\leq\infty$ we let $\mathrm{L}^p(G)$ be the usual L$^p$-space
with respect to $m$.  If $s\in G$ and $f$ is a complex-valued function on $G$
we denote the group action of left translation of $s$ on $f$ by $s\con f(t)
=f(s^{-1}t)$ for $t\iin G$.  For any appropriate pair of functions $f,g$ 
we denote $f\con g=\int_G f(s)s\con g \,ds$.

The {\it Fourier} and {\it Fourier-Stieltjes algebras},
$\falg$ and $\fsalg$ were defined by Eymard in \cite{eymard}. We recall
that $\falg$ consists of functions on $G$ of the form 
$s\mapsto\inprod{\lam(s)f}{g}=\bar{g}\con\check{f}(s)$, where
$\lam:G\to\bdop{\bltwog}$ is the left regular representation where
$\lam(s)f(t)=f(s^{-1}t)$ for almost every $t\iin G$.  We note that
$\falg$ has Gelfand spectrum $G$, given by evaluation.  We also
remark that we have duality relations $\falg^*\cong\vng$, where
$\vng$ is the von Neumann algebra generated by $\lam(G)$, and 
$\cstarg^*\cong\fsalg$ where $\cstarg$ is the universal C*-algebra
of $G$.

\subsection{Operator Spaces}\label{ssec:opspace}
Our main reference for operator spaces is \cite{effrosrB}.
An operator space is a complex vector space $\fV$, equipped with a sequence
of norms, one on each space of $n\cross n$ matrices over $\fV$, $\matn{\fV}$,
which satisfy Ruan's axioms; we call this an {\it operator space structure}.
An operator space is complete if $\mat{1}{\fV}$
is complete, i.e.\ $\fV=\mat{1}{\fV}$ is a Banach space.
If $\fW$ is another operator space and $T:\fV\to\fW$ is a linear map
we let $T^{(n)}:\matn{\fV}\to\matn{\fW}$ be the amplification
given by $T^{(n)}[v_{ij}]=[Tv_{ij}]$.  We say $T$ is {\it completely bounded}
if $\cbnorm{T}=\sup_{n\in\En}\norm{T^{(n)}}<\infty$.  Moreover we say
$T$ is a {\it complete contraction/isometry/quotient map} if each $T^{(n)}$
is a contraction/isometry/quotient map.  The space of completely bounded
linear maps from $\fV$ to $\fW$ is denoted $\cbop{\fV,\fW}$.

We make note of some basic operator space constructions.  If $\fW$ is a closed
subspace of an operator space $\fV$, then $\fW$ inherits the operator space
structure from $\fV$.  Moreover, the quotient space $\fV/\fW$ obtains the
{\it quotient} operator space structure via isometric identifications
$\matn{\fV/\fW}\cong\matn{\fV}/\matn{\fW}$, i.e.\ 
\[
\norm{[v_{ij}+\fW]}_{\matn{\fV/\fW}}=
\inf\{\norm{[v_{ij}+w_{ij}]}_{\matn{\fV}}:[w_{ij}]\in\matn{\fW}\}.
\]  
If $\fV$ and $\fW$ are any two
operator spaces, then the space $\cbop{\fV,\fW}$ obtains the {\it standard}
operator space structure (see \cite{blecher}) where we identify, for each $n$,
the matrix
$[S_{ij}]\iin\matn{\cbop{\fV,\fW}}$ with the operator $v\mapsto[S_{ij}v]$ in
$\cbop{\fV,\matn{\fW}}$.  We will make extensive use of the {\it operator
projective tensor product}, defined in \cite{effrosr}.  If $\fV$ and
$\fW$ are two complete operator spaces, let $\fV\optens\fW$ denote their
operator projective tensor product.  The algebraic tensor product
of $\fV$ and $\fW$ forms a dense subspace of $\fV\optens\fW$ and we let
$\fV\otimes_\wedge\fW$ denote
this algebraic tensor product, given the operator space
projective structure.

We note that for any operator space $\fV$ that every continuous 
linear functional is automatically completely bounded, i.e.\ 
$\fV^*=\cbop{\fV,\Cee}$, and thus is an operator space with the
standard operator space structure.  For a locally compact group $G$, 
the space $\fsalg\cong\cstarg^*$ is always endowed with the standard operator 
space structure.  $\falg$ obtains the same operator space structure as a 
subspace of $\fsalg$ as it does as the predual of $\vng$, i.e.\ a subspace of 
the dual.  $\bloneg$, as the predual of a commutative von Neumann algebra
naturally admits the {\it maximal} operator space structure.

Let $\fA$ be a Banach algebra, equipped with an operator space structure
such that $\mat{1}{\fA}=\fA$ isometrically, and let $\fV$ be a left 
$\fA$-module which is also an operator space.  For $a\iin\fA$ we let
$m_a:\fV\to\fV$ be the module action map, $M_av=a\mult v$.  
Then $\fV$ is called
a {\it completely bounded Banach $\fA$-module} under any of the following 
equivalent assumptions:

\smallskip
\noindent 
{\bf (i)}  $\{M_a:a\in\fA\}\subset\cbop{\fV}$ and $a\mapsto M_a:\fA\to
\cbop{\fV}$ is completely bounded,

\smallskip
\noindent {\bf (ii)} there is $C>0$ such that for any pair of matrices
$[a_{ij}]\iin\matn{\fA}$ and $[v_{kl}]\iin\mat{m}{\fV}$, we have
$\norm{[a_{ij}v_{kl}]}_{\mat{nm}{\fV}}\leq 
C\norm{[a_{ij}]}_{\matn{\fA}}\norm{[v_{kl}]}_{\mat{m}{\fV}}$, 
and

\smallskip\noindent 
{\bf (iii)} the map $\fA\otimes_\wedge\fV\to\fV$, given on elementary tensors
by $a\otimes v\mapsto a\mult v$ is completely bounded.

\smallskip
\noindent Similar definitions can be given with right and bimodules.
We say $\fV$ is a {\it completely contractive Banach $\fA$-module} if the maps
in (i) and (iii) above are complete contractions and in (ii) we can set $C=1$.
We call $\fA$ a {\it completely bounded (contractive) Banach algebra}, if it
is a completely bounded (contractive) module over itself.
We note that $\falg$, $\fsalg$ and $\bloneg$ are all completely contractive
Banach algebras with their standard operator space structures.

Operator Segal algebras were introduced in \cite{forrestsw}.  Let
$\fA$ be a completely contractive Banach algebra. An {\it operator
Segal algebra} in $\fA$ is a dense left ideal $\seg\fA$, equipped with a
complete operator space structure such that

\smallskip
\noindent {\bf (OSA1)} $\seg\fA$ is a completely bounded Banach
$\fA$-module, and

\smallskip
\noindent {\bf (OSA2)} the injection $\seg\fA\hookrightarrow\fA$ is completely
bounded.

\smallskip
\noindent We further call $\seg\fA$ a {\it contractive operator
Segal algebra} in $\fA$ if the maps above are complete contractions.
However, by uniformly scaling the matricial norms of $\seg\fA$ with
a fixed small enough constant, we may make any operator
Segal algebra a contractive one, and we will not insist on doing so in the 
sequel. We note that $\seg\fA$ itself is a completely bounded Banach algebra.

\bigskip
We finish this section by
outlining an approach to operator spaces and completely bounded maps
which is dual to the traditional one.  Let, for the remainder of the section,
$\fV$ and $\fW$ be complete operator spaces.

We let for any $n\iin\En$, $\stmatn$ denote the operator space of
$n\cross n$ matrices with the dual operator space structure, i.e., 
$\stmatn\cong\smatn^*$ completely isometrically.  
We let $\tmatn{\fV}=\stmatn\optens\fV$,
which we regard as matrices.  If 
$S:\fV\to\fW$ is a completely bounded map, then we let
$\tmatn{S}:\tmatn{\fV}\to\tmatn{\fW}$ denote the amplification, i.e.\
$\tmatn{S}=\mathrm{id}_{\stmatn}\otimes S$.
We also let $\stmatinfty$ denote the set of $\En\cross\En$ matrices
which may be identified as trace class operators on $\ell^2(\En)$,
which we endow with the usual predual operator space structure
$\stmatinfty\cong\fB(\ell^2(\En))_*$.  We define $\tmatinfty{\fV}$ analogously
as above, and also the operator $\tmatinfty{S}$, when it is defined.
We note the following elementary, but important fact.

\begin{tmatop}\label{prop:tmatop}
If $S:\fV\to\fW$ is a linear map, then the following are equivalent:

{\bf (i)} $S$ is a complete contraction [resp.\ complete quotient map],

{\bf (ii)} each $\tmatn{S}$ is a contraction [resp.\ quotient map], and

{\bf (iii)} $\tmatinfty{S}$ is defined and is a contraction [resp.\ quotient
map].
\end{tmatop}

\proof {\bf (i)} $\iff$ {\bf (ii)}.  
This is \cite[4.18]{effrosrB}, in light of the
identification \cite[(7.1.90)]{effrosrB} -- which shows that our
definition of $\tmatn{\fV}$ coincides with theirs.  
If $S$ is a complete quotient, then each $\tmatn{S}$ is a quotient map
by the projectivity property of the operator projective tensor product;
see \cite[7.1.7]{effrosrB}.  

{\bf (ii)} $\iff$ {\bf (iii)} for contractions.  
We have for each $n$ a completely isomorphic embedding
$\stmatn\hookrightarrow\stmatinfty$, given by identifying elements of
$\stmatn$ with elements of $\stmatinfty$ whose non-zero entries are only
in the upper left $n\cross n$ corner.  Since 
$\stmatinfty^*\cong\bdop{\ell^2(\En)}$
is an injective operator space, we obtain completely isometric imbeddings
$\tmatn{\fV}=\stmatn\optens\fV\hookrightarrow\tmatinfty{\fV}
=\stmatinfty\optens\fV$;
see the discussion \cite[p.\ 130]{effrosrB}.  As $\bigcup_{n=1}^\infty\stmatn$
is dense in $\stmatinfty$, it follows that 
$\tmatfin{\fV}=\bigcup_{n=1}^\infty\tmatn{\fV}$
is dense in $\tmatinfty{\fV}$.  Now we can define $\tmatfin{S}:\tmatfin{\fV}\to
\tmatfin{\fW}$ in the obvious way -- so $\tmatfin{S}=
\tmatinfty{S}|_{\tmatinfty{\fV}}$, when the latter makes sense.  
We have that
\[
\norm{\tmatfin{S}}=\sup_{n\in\En}\norm{\tmatn{S}}.
\]
Thus if (ii) is assumed, then $\tmatfin{S}$ is contractive, and thus
$\tmatinfty{S}$ is defined and contractive.  Conversely, if (iii) is assumed
than $\tmatfin{S}$ is contractive,whence (ii).  

{\bf (ii)} $\iff$ {\bf (iii)} for quotient maps.
Suppose that $\tmatinfty{S}$ is a quotient map.
By \cite[10.1.4]{effrosrB}, we have that $\tmatinfty{\fV}^*\cong
\mat{\infty}{\fV^*}$, where
$\mat{\infty}{\fV^*}$ is the space of $\En\cross\En$ matrices with entries
in $\fV^*$ whose finite submatrices are uniformly bounded in norm.
Thus $\tmatinfty{S}^*=(S^*)^{(\infty)}:\mat{\infty}{\fW^*}\to
\mat{\infty}{\fV^*}$
is an isometry.  It follows that $S^*$ is a complete isometry, whence 
$S$ is a complete quotient map by \cite[4.1.8]{effrosrB}. \endpf

We say that a linear operator $S:\fV\to\fW$ is a {\it complete isomorphism}
if it is completely bounded, bijective, and $S^{-1}:\fW\to\fV$ is 
completely bounded too.  We say that $S:\fV\to\fW$ is a {\it complete 
surjection} if the induced map $\til{S}:\fV/\ker S\to\fW$,
defined by $\til{S}q=S$ where $q:\fV\to\fV/\ker S$ is the quotient map,
is a complete isomorphism.

\begin{tmatop1}\label{cor:tmatop1}
{\bf (i)} Suppose $S\iin\cbop{\fV,\fW}$ is a bijection.  
Then $S$ is a complete isomorphism
if and only if $\tmatinfty{S}:\tmatinfty{\fV}\to\tmatinfty{\fW}$ is 
an isomorphism of Banach spaces.  

{\bf (ii)} Suppose $S\iin\cbop{\fV,\fW}$ is surjective.  Then
$S$ is a complete surjection 
if and only if $\tmatinfty{S}:\tmatinfty{\fV}\to\tmatinfty{\fW}$ is 
surjective.\end{tmatop1}

\proof {\bf (i)} If $S^{-1}\in\cbop{\fW,\fV}$, then 
$\tmatinfty{S^{-1}}=\tmatinfty{S}^{-1}$.
Conversely, if $T=\tmatinfty{S}^{-1}$ is a bounded operator, we have that
for $w\iin\tmatfin{\fW}$ that
\[
\tmatinfty{S}Tw=w=\tmatinfty{S}\tmatfin{S^{-1}}w
\]
so $T|_{\tmatfin{\fW}}=\tmatfin{S^{-1}}|_{\tmatfin{\fW}}$.  Thus
$\tmatfin{S^{-1}}$ is bounded, whence so too is $\tmatinfty{S^{-1}}$.  

{\bf (ii)}  If $S$ is surjective than $\til{S}$ is bijective.  From above,
if $\til{S}$ is a complete isomorphism, then
$\tmatinfty{\til{S}}$ is an isomorphism of Banach spaces.  It follows
that $\tmatinfty{S}=\tmatinfty{\til{S}}\tmatinfty{q}$ is surjective.
On the other hand, if $\tmatinfty{S}=\tmatinfty{\til{S}}\tmatinfty{q}$ is 
surjective
then $\tmatinfty{\til{S}}$ is surjective.  As it is already injective, 
and bounded as $\tmatinfty{q}$ is a quotient map, we obtain that
$\tmatinfty{\til{S}}$ is a bounded bijection, hence an isomorphism
by the open mapping theorem.   \endpf

We have, by \cite[7.1.6]{effrosrB} that
$\stmatn\optens\fW^*=\tmatn{\fW^*}\cong\matn{\fW}^*\cong(\smatn\stens\fW)*$,
where the duality is given in tensor form by $\dpair{t\otimes f}{m\otimes w}
=\mathrm{trace}(tm)f(w)$, for $t\in\stmatn$, $f\in\fW^*$, $m\in\smatn$ and
$w\in\fW$.  Thus, in matrix form, this dual paring becomes
\[
\bigl\langle [f_{ij}],[w_{ij}]\bigr\rangle=\sum_{i,j=1}^nf_{ij}(w_{ji})
\]
for $[f_{ij}]\iin\tmatn{\fW^*}$ and $[w_{ij}]\iin\matn{\fW}$.  Thus the map
$[S_{ij}]$ in $\cbop{\fV,\matn{\fW}}$ has adjoint 
$[S_{ij}]^*\iin\cbop{\tmatn{\fW^*},\fV^*}$
given by
\[
\bigl\langle [S_{ij}]^*[f_{ij}],v\bigr\rangle=
\bigl\langle [f_{ij}],[S_{ij}v]\bigr\rangle=\sum_{i,j=1}^nf_{ij}(S_{ji}v)
=\sum_{i,j=1}^nS_{ji}^*f_{ij}(v)
\]
for $[f_{ij}]\iin\tmatn{\fW^*}$ and $v\in\fV$. 

Now if $[S_{ij}]\in\cbop{\fV,\matn{\fW}}$,
we have that $[S_{ij}^*]\in\cbop{\fW^*,\matn{\fV^*}}$, with 
$\cbnorm{[S_{ij}^*]}
=\cbnorm{[S_{ij}]}$, 
with proof similar to that of \cite[3.1.2]{effrosrB}.
Then $[S_{ij}^*]^*\in\cbop{\tmatn{\fV^{**}},\fW^{**}}$ and we let
\begin{equation}\label{eq:sss}
[S_{ij}^*]_*=[S_{ij}^*]^*|_{\tmatn{\fV}}\in\cbop{\tmatn{\fV},\fW}.
\end{equation}
Thus if $[v_{ij}]\in\tmatn{\fV}$, we have obtain 
$[S_{ij}^*]_*[v_{ij}]=\sum_{i,j=1}^nS_{ji}v_{ij}$.
We have adjoint ${[S_{ij}^*]_*}^*=[S_{ij}^*]$, and hence 
\begin{equation}\label{eq:sssnorm}
\cbnorm{[S_{ij}^*]_*}=\cbnorm{[S_{ij}^*]}=\cbnorm{[S_{ij}]}.
\end{equation}
This equation will be useful in the sequel when we determine that
$\sofalg$ is an operator Segal algebra in $\falg$.

\subsection{Localisation}\label{ssec:localisation}
Let $\fA$ be a semi-simple commutative Banach algebra with Gelfand spectrum
$X$.  Via the Gelfand transform, we regard $\fA$ as a subalgebra
of $\mathcal{C}_0(X)$, the algebra of continuous functions on $X$
vanishing at infinity.  We say that $\fA$ is
{\it regular} if for every pair $(x,F)$, where 
$x\in X$ and $F$ is a closed subset of $X$
with $x\not\in F$,
we have that there is $u\iin\fA$ such that $u(x)=1$ and $u|_F=0$.
We note, below, that such an algebra admits local inverses.

\begin{standard}
If $\fA$ is a regular Banach algebra with Gelfand spectrum $X$,
$u\in\fA$ and $K$ is a compact subset of $X$ on which $u$ does not vanish,
then there is $u'\iin\fA$ such that $uu'|_K=1$.
\end{standard}

\proof By \cite[(39.12)]{hewittrII}, $\fA|_K=\{u|_K:u\in\fA\}$
is a regular Banach algebra with Gelfand spectrum $K$.
Then we may apply analytic functional calculus (\cite[39.14]{hewittrII})
to see that $1/(u|_K)\in\fA|_K$.  Find $v\iin\fA$ such that $v|_K=1/(u|_K)$.
\endpf

Thus, we see that a regular Banach algebra $\fA$ is a ``standard function
algebra on its spectrum'', in the sense of Reiter~\cite{reiter,reiters}.
Let $\fI$ be an ideal in $\fA$, which is not necessarily assumed to be closed.
We define the {\it hull} of $\fI$ to be the set $\mathrm{h}(\fI)=
\{x\in X:u(x)=0\text{ for each }u\in\fI\}$.
We let $\fA_c=\{u\in\fA:\supp{u}\text{ is compact}\}$.  
Thus we obtain the following localisation result, 
proved in \cite[2.1.4]{reiter} (or \cite[2.1.14]{reiters}),
which we restate here, without proof, for convenient reference.

\begin{localisation}\label{cor:localisation}
Let $\fA$ be a regular Banach algebra and $\fI$ be an ideal of $\fA$.  
If $u\in\fA_c$ with $\supp{u}\cap\mathrm{h}(\fI)=\varnothing$, then
$u\in\fI$.  In particular, if $\mathrm{h}(\fI)=\varnothing$, then
$\fA_c\subset\fI$.
\end{localisation}

We say that an ideal $\fI$ of a regular
Banach algebra $\fA$ has {\it compact support}, if $\supp{\fI}
=\wbar{G\setdif\mathrm{h}(\fI)}$ is compact in the spectrum $X$.
The following result can be applied
in more general situations than we give, but is only required for $\falg$
where $G$ is a locally compact group.

\begin{localisation1}\label{cor:localisation1}
Let $\fI$ and $\fJ$ each be non-zero ideals having compact support in
$\falg$ with $\fI\subset\fJ$.  Then there are $t_1,\dots,t_n\iin G$
and $u_1,\dots,u_n\iin\fI$ such that
\[
\sum_{l=1}^n(t_l\con u_l)v=v\text{ for all }v\iin\fJ.
\]
\end{localisation1}

\proof Let $Q$ be any compact subset of $\supp{\fI}$ having
non-empty interior $Q^\circ$.  Find $t_1,\dots,t_n\iin G$ so
$\bigcup_{l=1}^nt_lQ^\circ\supset\supp{\fJ}$.  Then
$\fI'=\sum_{l=1}^nt_l\con\fI$ is an ideal in $\falg$
for which $\supp\fJ\subset(\supp\fI)^\circ$.  
By the the regularity of $\falg$ and \cite[(39.15)]{hewittrII}
(or see \cite[(3.2)]{eymard}), there is a function $u\in\falg$
such that $u|_{\supp{\fJ}}=1$ and $\supp{u}\subset\bigcup_{l=1}^nt_lQ$.
By Corollary \ref{cor:localisation}, above, $u\in\fI'$, hence
$u=\sum_{l=1}^nt_l\con u_l$, as desired.  \endpf

\section{Construction of Feichtinger's Segal Algebra}

Let $G$ be a locally compact group.
In this section, we reconstruct Feichtinger's
Segal algebra, which was done explicitly for
abelian $G$ in \cite{feichtinger} and \cite[\S 6.2]{reiters} and, implicitly
for general $G$, in \cite{feichtinger81}.
The construction we give below is superficially
different than the one of Feichtinger, but conceptually
more useful to our task.

Let $\fI$ be a non-zero closed ideal in $\falg$ which has compact
support, as defined in Section \ref{ssec:localisation}.
We let $\ell^1(G)$ be the usual
$\ell^1$-space, indexed over $G$, which will be identified
with the closed linear span of the Dirac measures $\{\del_s:s\in G\}$.
We provide $\ell^1(G)$ with the maximal operator space structure.
We define
\[
q_\fI:\ell^1(G)\optens\fI\to\falg\quad
\text{ by }\quad q_\fI(x\otimes u)=x\con u=\sum_{s\in G}\alp_ss\con u
\]
where $x=\sum_{s\in G}\alp_s\del_s$ and
$s\con u(t)=u(s^{-1}t)$ for $t\iin G$.  
We let 
\[
\sofalg=\ran q_\fI\cong\ell^1(G)\optens\fI/\ker q_\fI.
\]
We make $\sofalg$ into an operator space by giving it
the quotient operator space structure. 
Let us note the following description of $\sofalg$, which is
exactly that of \cite{feichtinger}, when $n=1$ and
$\fI=\falkg$, where
\[
\falkg=\{u\in\falg:\supp{u}\subset K\}
\]
for $K$ a compact subset of $G$ with non-empty interior.
The matrices $\tmatn{\fV}$, for an operator space $\fV$, are
described in Section \ref{ssec:opspace}, above.

\begin{feichdesc}\label{lem:feichdesc}
We have for any $n=1,2,\dots,\infty$ that
\begin{equation}\label{eq:feichdesc}
\tmatn{\sofalg}=
\left\{\sum_{k=1}^\infty [s_k\con u_{ij}^{(k)}]:
\begin{matrix} s_k\in G,[u_{ij}^{(k)}]\in\tmatn{\fI}\wwith \\
\sum_{i=1}^\infty\norm{[u_{ij}^{(k)}]}_{\tmatn{\mathrm{A}}}<+\infty
\end{matrix}\right\}.
\end{equation}
The norm on $\tmatn{\sofalg}=\tmatn{\ran q_\fI}$ is given by
\[
\norm{[u_{ij}]}_{\tmatn{\ran q_\fI}}=
\inf\left\{\sum_{k=1}^\infty\norm{[u_{ij}^{(k)}]}_{\tmatn{\mathrm{A}}}:
[u_{ij}]=\sum_{k=1}^\infty [s_k\con u_{ij}^{(k)}]\text{ as above}\right\}.
\]
\end{feichdesc}

\proof First, let $n=1$.
We observe that there is an isometric identification
$\ell^1(G)\optens\fI=
\ell^1(G)\ptens\fI$, by virtue of the fact that $\ell^1(G)$ has maximal 
operator space structure; see \cite[8.2.4]{effrosrB}.  Now
if $t\in\ell^1(G)\ptens\fI$ and $\eps>0$, we can write
\[
t=\sum_{k=1}^\infty\left(\sum_{s\in G}\alp_{ks}\del_s\right)\otimes u^{(k)},
\wwhere\sum_{k=1}^\infty\left(\sum_{s\in G}|\alp_{ks}|\right)
\fnorm{u^{(k)}}<\norm{t}_\gam+\eps.
\]
We see the sum can easily be rearranged in the form
$t=\sum_{s\in G}\del_s\otimes\left(\sum_{k=1}^\infty\alp_{ks}u^{(k)}\right)$.
(This is essentially the proof that $\ell^1(G)\ptens\fI$ is isometrically
isomorphic to the $\fI$-valued $\ell^1$-space, $\ell^1(G;\fI)$.)
Thus we see that equations for $u\iin\sofalg$ can be arranged as
suggested in (\ref{eq:feichdesc}), and any such equation describes
an element of $\sofalg=\ran q_\fI$.
The formula for the norm follows immediately from the fact that
$\sofalg$ is a quotient of $\ell^1(G)\otimes^\gam\fI$.  

Now if $n>1$, we make identifications
\[
\tmatn{\ell^1(G)\optens\fI}\cong\stmatn\optens\bigl(\ell^1(G)\optens\fI\bigr)
\cong\ell^1(G)\optens\bigl(\stmatn\optens\fI\bigr)
=\ell^1(G)\ptens\tmatn{\fI}.
\]
The proof then follows the $n=1$ case, given above.  \endpf

We now state some of the remarkable properties of $\sofalg$.

\begin{maintheorem}\label{theo:maintheorem}
{\bf (i)} The space $\sofalg$ is a (contractive)
operator Segal algebra in $\falg$.

{\bf (ii)} $\sofalg$ is the smallest Segal algebra $\sfalg$ in $\falg$  
which is closed under left translations and on which left translations
are isometric, i.e.\ $s\con u\in\sfalg$ for each $s\iin G$ and $u\iin\sfalg$
with $\sfnorm{s\con u}=\sfnorm{u}$.  Moreover, for each
$u\in\sofalg$, $s\mapsto s\con u:G\to\sofalg$ is continuous, and for each
$s\iin G$, $u\mapsto s\con u:\sofalg\to\sofalg$ is a complete isometry.

{\bf (iii)} For any closed ideals $\fI,\fJ$ of $\falg$, each having compact 
support, the operator space structure on $\sofalg$ qua $\ran q_\fI$, or on 
$\sofalg$ qua $\ran q_\fJ$, are completely isomorphic.
\end{maintheorem}

\proof {\bf (i)}  As in the construction above, we fix a non-zero closed
ideal in $\falg$ with compact support.

We note that it is easy to see, using 
Lemma \ref{lem:feichdesc} above, and Proposition \ref{prop:tmatop},
that $\sofalg$ imbeds completely contractively into $\falg$. 
Moreover, using the $n=1$ case of the lemma, it is easy to see that
$\sofalg$ is a Banach $\falg$-module. Unfortunately, it is somewhat involved 
to show that $\sofalg$ is a completely contractive $\falg$-module.

Let for $v\iin\falg$, $M_v:\falg\to\falg$ be the multiplication map.
Using (\ref{eq:sss}) and (\ref{eq:sssnorm}), it suffices 
to show that for any $[v_{ij}]\in\matn{\falg}$ we have that
\[
[M_{v_{ij}}^*]_*\in\cbop{\tmatn{\sofalg},\sofalg}
\wwith\cbnorm{[M_{v_{ij}}^*]_*}\leq\norm{[v_{ij}]}_{\matn{\mathrm{A}}}.
\]
By Proposition \ref{prop:tmatop} it suffices to show for each $m$
that
\[
\tmatm{[M_{v_{ij}}^*]_*}:\tmatm{\tmatn{\sofalg}}\to\tmatm{\sofalg}
\]
satisfies 
\begin{equation}\label{eq:strangenorm}
\norm{\tmatm{[M_{v_{ij}}^*]_*}}\leq
\norm{[v_{ij}]}_{\matn{\mathrm{A}}}.
\end{equation}
We note that there is a natural identification
\[
\tmatm{\tmatn{\sofalg}}\cong\tmatnm{\sofalg}
\]
of which we take advantage.  However, we will prefer,
for computational convenience, to label elements of this space by doubly
indexed matrices $[u_{ij,pq}]$, where $i,j=1,\dots,n$ and $p,q=1,\dots,m$.

By Lemma \ref{lem:feichdesc}, each element $[u_{ij,pq}]$ of
$\tmatnm{\sofalg}$ admits, for any $\eps>0$, the form
\[
[u_{ij,pq}]=\sum_{k=1}^\infty [s_k\con u_{ij,pq}^{(k)}]=\tmatn{q_\fI}
\left(\sum_{k=1}^\infty
\del_{s_k}\otimes[u_{ij,pq}^{(k)}]\right)
\]
where 
$\sum_{k=1}^\infty\norm{[u_{ij,pq}^{(k)}]}_{\tmatnm{\mathrm{A}}}
\leq\norm{[u_{ij}]}_{\tmatnm{\ran q_\fI}}+\eps$.
We see that
\begin{align*}
\tmatm{[M_{v_{ij}}^*]_*}[w_{ij,pq}]
&=\sum_{k=1}^\infty \tmatm{[M_{v_{ij}}^*]_*}
\bigl([s_k\con u_{ij,pq}^{(k)}]\bigr) \\
&=\sum_{k=1}^\infty
\left[s_k\con\sum_{i,j=1}^m (s_k^{-1}\con v_{ji})u_{ij,pq}\right] \\
&=\tmatm{q_\fI}
\left(\sum_{k=1}^\infty\del_{s_k}\otimes
\left[\sum_{i,j=1}^m (s_k^{-1}\con v_{ji})u_{ij,pq}\right]\right) \\
&=\tmatm{q_\fI}\left(\sum_{k=1}^\infty\del_{s_k}\otimes
\tmatm{[M_{s_k^{-1}\ast v_{ij}}^*]_*}[u_{ij,pq}^{(k)}]\right).
\end{align*}
Using Proposition \ref{prop:tmatop}, (\ref{eq:sssnorm}), and the fact
that translation is a complete isometry on $\falg$, we thus obtain
\begin{align*}
&\norm{\sum_{k=1}^\infty\del_{s_k}\otimes\tmatm{[M_{s_k^{-1}\ast v_{ij}}^*]_*}
[u_{ij,pq}^{(k)}]}_{\tmatnm{\ell^1\optens\fI}} \\
&\qquad\leq\sum_{k=1}^\infty\norm{ \tmatm{ [M_{ s_k^{-1}\ast v_{ij} }^*]_* } }
\norm{[u_{ij,pq}^{(k)}]}_{\tmatnm{\mathrm{A}}} \\
&\qquad=\sum_{k=1}^\infty\norm{ [s_k^{-1}\con v_{ij}] }_{\matm{\mathrm{A}}}
\norm{[u_{ij,pq}^{(k)}]}_{\tmatnm{\mathrm{A}}} \\
&\qquad\leq\norm{[v_{ij}]}_{\matn{\mathrm{A}}}
\bigl(\norm{[u_{ij,pq}]}_{\tmatnm{\ran q_\fI}}+\eps\bigr).
\end{align*}
Thus, since $\eps$ is arbitrary, we obtain that 
\[
\norm{\tmatm{[M_{v_{ij}}^*]_*}[u_{ij,pq}]}_{\tmatm{\ran q_\fI}}\leq
\norm{[v_{ij}]}_{\matn{\mathrm{A}}}\norm{[u_{ij,pq}]}_{\tmatnm{\ran q_\fI}}
\]
and hence we obtain (\ref{eq:strangenorm}).

{\bf (ii)} It it immediate from the construction of $\sofalg$ that
it is closed under left translations.
We note that the action of $G$ on $\sofalg$ is continuous, and
is one of isometries, in fact complete isometries.  This follows by a 
straightforward
application of Lemma \ref{lem:feichdesc} and Proposition \ref{prop:tmatop}.

Any Segal algebra $\sfalg$ in $\falg$ is an ideal with
empty hull, and thus, by Corollary \ref{cor:localisation},
necessarily contains $\falcg$.  Hence for any non-zero compactly 
supported ideal
$\fI$ of $\falg$, we obtain that
$\fI\subset\sfalg$.  Thus, if translations are pointwise continuous and 
isometric on $\sfalg$, we see by Lemma \ref{lem:feichdesc}, in the case $n=1$,
that $\sfalg\supset\sofalg$.

{\bf (iii)}  It follows from (ii), above, that 
$\sofalg$ is independent if the choice of ideal $\fI$.
Let us suppose that $\fI$ and $\fJ$ are two closed non-zero
ideals of $\falg$ having compact supports.  By replacing
$\fI$ by $\fI\cap\fJ$, if necessary, we may suppose $\fI\subset\fJ$.
Then the injection $\iota:\ell^1(G)\optens\fI\hookrightarrow
\ell^1(G)\optens\fJ$ is a complete contraction.
It is clear that $q_\fJ\comp\iota= q_\fI$, so $\iota$ induces
a completely contractive map $\til{\iota}:\ran q_\fI\to\ran q_\fJ$, which is
the identity map. 
Thus, by Proposition \ref{prop:tmatop}, 
$\tmatinfty{\til{\iota}}:\tmatinfty{\ran q_\fI}\to
\tmatinfty{\ran q_\fJ}$ is a contraction.  Let us see that 
$\tmatinfty{\til{\iota}}$ is 
surjective.  Let $u=\sum_{l=1}^nt_l\con u_l$ be as in Corollary 
\ref{cor:localisation1}.
We note that if $[w_{ij}]\in\tmatinfty{\ran q_\fJ}$, has form 
$[w_{ij}]=\sum_{k=1}^\infty[s_k\con w_{ij}^{(k)}]$ as in (\ref{eq:feichdesc}),
 then we have
\begin{align}\label{eq:surjectivity}
[w_{ij}]&=\sum_{k=1}^\infty[s_k\con (uw_{ij}^{(k)})] 
=\sum_{k=1}^\infty\sum_{l=1}^n[s_k\con (t_l\con u_l\,w_{ij}^{(k)})] \\ \notag
&=\tmatinfty{q_\fI}\left(\sum_{k=1}^\infty\sum_{l=1}^n\del_{s_kt_l}\otimes
[u_l(t_l^{-1}\con w_{ij}^{(k)})]\right)  
\end{align}
which is an element of $\tmatinfty{\ran q_\fI}$.
Hence $\tmatinfty{\til{\iota}}$ is surjective, and thus,  
by the open mapping theorem,
an isomorphism of Banach spaces.  We then appeal to Corollary 
\ref{cor:tmatop1}.
\endpf

Let us note that Theorem \ref{theo:maintheorem}, above, holds under more 
general assumptions.  Though the assumptions we give below seem
less natural, they are indispensable for actually working with $\sofalg$.
Let $\fI$ be a non-zero compactly supported ideal in $\falg$, 
which is not necessarily closed, but comes equipped
with an operator space structure by which it is a completely contractive
Banach $\falg$-module, and the inclusion map $\fI\hookrightarrow\falg$
is completely bounded.  
We define $q_\fI:\ell^1(G)\optens\fI\to\falg$, and
its quotient space $\ran q_\fI$, with its quotient operator space structure, 
as before.

\begin{maintheorem1}\label{cor:maintheorem1}
With $\fI$ as above, $\ran q_\fI$ is an operator Segal algebra in $\falg$
which is completely isomorphic with
$\ran q_{\bar{\fI}}$, where $\bar{\fI}$ is the closure of $\fI$ in $\falg$.
Hence $\ran q_\fI=\sofalg$, completely isomorphically.
\end{maintheorem1}

\proof We first note that the proofs of Lemma \ref{lem:feichdesc},
and then of (i) and (ii) of the theorem above,
can be applied verbatim; though we should note that the inclusion
$\ran q_\fI\hookrightarrow\falg$ is completely bounded, instead of
completely contractive.  Thus we see that $\ran q_\fI$ is an 
operator Segal algebra in $\sofalg$ on which $G$ acts
continuously and isomorphically by translations.  Hence
$\ran q_\fI=\ran q_{\wbar{\fI}}$.  Moreover
the proof of part (iii) can be applied up to seeing that 
$\tmatinfty{\til{\iota}}:\tmatinfty{\ran q_\fI}\to
\tmatinfty{\ran q_{\bar{\fI}}}$ is   
bounded.  To see obtain surjectivity, we note for any $u\iin\fI$
that $M_u:\falg\to\fI$ is completely bounded, and then (\ref{eq:surjectivity}),
where $[w_{ij}]\in\tmatinfty{\ran q_{\bar{\fI}}}$, has the form
\[
[w_{ij}]=\tmatinfty{q_\fI}\left(\sum_{k=1}^\infty\sum_{l=1}^n\del_{s_kt_l}
\otimes
\tmatinfty{M_{u_l}}[t_l^{-1}\con w_{ij}^{(k)}]\right) 
\]
where each $\tmatinfty{M_{u_l}}[t_l^{-1}\con w_{ij}^{(k)}]\in\tmatinfty{\fI}$.
\endpf

We recall that a Segal algebra $\soneg$ in $\bloneg$ is called {\it symmetric}
if for any $s\iin G$ and $f\iin\soneg$ we have
\[
\sonenorm{s\con f}=\sonenorm{f}=\sonenorm{f\con s}
\]
where $f\con s(t)=\Del(s)^{-1}f(ts^{-1})$ for almost every $t\iin G$.  This
is a necessary and sufficient condition to make $\soneg$ a two-sided
ideal in $\bloneg$.  $\soneg$ is called {\it pseudo-symmetric} if the
anti-action $s\mapsto f\con s$ is continuous on $G$ for any fixed 
$f\iin\soneg$.
We note that this is equivalent to having the action of right translation
\begin{equation}\label{eq:rtrans}
s\mapsto s\mult f,\wwhere s\mult f(t)=f(ts)\text{ for almost every }t 
\end{equation}
continuous on $G$ for any fixed $f\iin\soneg$.
We also recall, by well known technique (see \cite[pps.\ 26-27]{johnsonm}, 
for example), that any Banach space $\fV$ is an contractive essential 
left/right $\bloneg$-module if and only if there
is a continuous action/anti-action of $G$ on $\fV$ by linear isometries.
Furthermore, if $\fV$ is an operator space, $\fV$ is a completely
contractive $\bloneg$-module, if and only if the associated action
of $G$ on $\fV$ is one by complete isometries.

It will be useful, below, to recall the {\it Lebesgue-Fourier algebra}
\begin{equation}\label{eq:lebfalg}
\lebfalg=\falg\cap\bloneg
\end{equation}
studied in \cite{ghahramanil,ghahramanil2,forrestsw}.  This is simultaneously
a Segal algebra in $\falg$ and in $\bloneg$.  In \cite{forrestsw} it was shown
that $\lebfalg$ is also a contractive operator Segal algebra in either context.

\begin{maintheorem2}\label{cor:maintheorem2}
{\bf (i)} $\sofalg$ is a pseudo-symmetric operator
Segal algebra in $\bloneg$.  It is symmetric only if $G$
is unimodular.

{\bf (ii)} Let $\fI$ be a fixed ideal in $\falg$ satisfying the assumptions of
Corollary \ref{cor:maintheorem1}.  Let $q_\fI':\bloneg\optens\fI\to\falg$
be given, on elementary tensors, by $q_\fI'(f\otimes u)=f\con u$.
Then $\ran q_\fI'$, with its quotient operator space structure, is
completely isomorphic to $\sofalg$.
\end{maintheorem2}

\proof {\bf (i)} If $\fI$ is a closed compactly supported pointwise
ideal in $\lebfalg$.  Then $\fI$ imbeds completely contractively into
$\bloneg$, and is a completely contractive $\falg$-module.  
Thus, it follows Corollary \ref{cor:maintheorem1} above, that
$\ran q_\fI$, with its quotient operator space structure, imbeds
completely contractively into $\bloneg$.  It then follows part (ii) of Theorem
\ref{theo:maintheorem} that $\sofalg$ is an operator Segal algebra
in $\bloneg$; and the same proof of can be trivially adapted to see $G$
acts continuously and completely isometrically on $\sofalg$
by right translation.  It is clear that $s\mult u=\Del(s)u\con s^{-1}$ for all
$u\iin\sofalg$ and $s\iin G$, and hence $\sofalg$ is a symmetric
Segal algebra in $\bloneg$, if and only if $G$ is unimodular.

{\bf (ii)} Let $(e_U)$ be the bounded approximate identity for $\bloneg$
given by normalised indicator functions of relatively compact
neighbourhoods of the identity, $e$.
If $[u_{ij}]\in\tmatinfty{\fI}$, then for $s\iin G$ we have that 
\[
[s\con u_{ij}]
=\lim_{U\searrow e}[s\con e_U\con u_{ij}]\in\tmatinfty{\ran q_\fI'}.
\]
Thus we see from Lemma \ref{lem:feichdesc}, that $\sofalg\subset
\ran q_\fI'$, completely boundedly.  
To obtain the converse inclusion, we note, similarly as in the proof
of Lemma \ref{lem:feichdesc}, that
\[
\tmatinfty{\bloneg\optens\fI}\cong\bloneg\optens\tmatinfty{\fI}
=\bloneg\ptens\tmatinfty{\fI}.
\]
Hence every element of $\ran q_\fI'$ is of the form
\begin{equation}\label{eq:feichdesc1}
\sum_{k=1}^\infty [f_k\con u_{ij}^{(k)}]
\wwhere \sum_{k=1}^\infty\blonenorm{f_k}
\norm{[u_{ij}^{(k)}]}_{\tmatinfty{\mathrm{A}}}<+\infty.
\end{equation}
Now if $[u_{ij}]\in\tmatinfty{\fI}$ and $f\in\bloneg$, then
\[
[f\con u_{ij}]=\int_G f(s)[s\con u_{ij}]\in\tmatinfty{\sofalg}
\]
since the integral may be realised as a Bochner integral in
$\tmatinfty{\sofalg}$, by the continuity of 
$s\mapsto[s\con u_{ij}]:G\to\tmatinfty{\sofalg}$.
Thus (\ref{eq:feichdesc1}) shows that $\ran q_\fI'\subset\sofalg$,
completely boundedly. \endpf

The next result shows the only occasions for which we know that
$\sofalg=\lebfalg$, and we conjecture these are all such occasions.
For this result, we will consider $\sofalg$ as a Segal algebra
in $\falg$.

\begin{maintheorem3}\label{cor:maintheorem3}
If $K$ is an open compact subgroup, with $T$ a transversal for left cosets,
then there is a natural completely isomorphic algebra homomorphism
\[
\sofalg\cong\ell^1(T)\optens\fal{K}
\]
where $\ell^1(T)$ has pointwise multiplication.  In particular,

{\bf (i)} $\sofalg=\ell^1(G)$, completely isomorphically, if $G$ is discrete,
and

{\bf (ii)} $\sofalg=\falg$, completely isomorphically, if $G$ is compact.
\end{maintheorem3}

\proof We first note that $\fal{K}\cong\falkg$, completely isometrically.
Second, for $s,t\iin G$, $s\con\falkg=\falskg$ and $t\con\falkg=\falbg{tK}$
are either identical or disjoint, depending on whether $s^{-1}t\in K$ or not.
Thus we can see that
$q_{\falkg}|_{\ell^1(T)\optens\falkg}:\ell^1(T)\optens\falkg\to
\ran q_{\falkg}$ is a complete isomorphism.  Indeed, it follows 
from Lemma \ref{lem:feichdesc} that 
\[
\tmatinfty{q_{\falkg}|_{\ell^1(T)\optens\falkg}} :
\tmatinfty{\ell^1(T)\optens\fal{K}}\to
\tmatinfty{\sofalg}
\]
is an bijection, hence an isomorphism.
\endpf

It is useful to observe the more general fact below.

\begin{maintheorem4}\label{cor:maintheorem4}
If $H$ is an open subgroup of $G$, with $T$ a transversal for left cosets,
then there is a natural completely isomorphic algebra homomorphism
\[
\sofalg\cong\ell^1(T)\optens\sofalh
\]
where $\ell^1(T)$ has pointwise multiplication.
\end{maintheorem4}

\proof Let $\fI$ be a non-empty closed compactly supported ideal of $\falg$
for which $\supp{\fI}\subset H$.  Then we may consider $\fI$ to be an ideal
in $\falh\cong\falbg{H}$, and we have that 
$q^H_\fI=q_\fI|_{\ell^1(H)\optens\fI}:\ell^1(H)\optens\fI\to\sofalh$ 
is a complete surjection.
The bijection $(t,s)\mapsto ts:T\cross H\to G$ induces an
isomorphism
\[
\ell^1(T)\optens\ell^1(H)=\ell^1(T)\ptens\ell^1(H)\cong
\ell^1(T\cross H)\cong\ell^1(G).
\]
Then we obtain the following commuting diagram.
\[
\xymatrix{
\ell^1(T)\optens\ell^1(H)\optens\fI\ar[rr]^{\cong}
\ar[d]_{\id\otimes q^H_\fI}
&  & \ell^1(G)\optens\fI\ar[d]^{q_\fI} \\
\ell^1(T)\optens\sofalh \ar[rr]^{\del_t\otimes u\mapsto t\ast u} & & \sofalg}
\]
We obtain, as in the proof of the result above, that 
the bottom arrow represents a complete isomorphism.  \endpf

\section{Functorial Properties}

\subsection{Tensor products}
Let us first note the primary motivation for desiring an
operator space structure on $\sofalg$.  This is an analogue
of a result from \cite{effrosr} which states that
\[
\falg\optens\falh\cong\falgh
\]
completely isometrically, via the natural morphism which identifies
$u\otimes v$ with the function $(s,t)\mapsto u(s)v(t)$.  Alternatively
we may view this as an analogue of the classical result that
\[
\bloneg\optens\blone{H}=\bloneg\ptens\blone{H}\cong\blone{G\cross H}
\]
where the projective and operator projective tensor products agree
since $\bloneg$ (or $\blone{H}$) is a maximal operator space.

\begin{tensorprod}\label{theo:tensorprod}
Let $G$ and $H$ be locally compact groups.
Then there is a natural complete isomorphism $\sofalg\optens\sofalh
\cong\sofalgh$.
\end{tensorprod}


\proof By \cite[(7.3)\&(7.4)]{hewittrI}, there are almost connected open
subgroups $G_0$ of $G$, and $H_0$ of $H$.  Let $\fI$ and $\fJ$ be
compactly supported ideals of $\fal{G_0}$ and $\fal{H_0}$, respectively.
The dual spaces $\fal{G_0}^*\cong\mathrm{VN}(G_0)$ and
$\fal{H_0}^*\cong\mathrm{VN}(H_0)$ are injective von Neumann algebras and
hence injective operator spaces; see \cite[pps.\ 227-228]{paterson}, for
example.  Thus, by comments in \cite[p.\ 130]{effrosrB}, the inclusion
maps induce a complete isometry
$\fI\optens\fJ\hookrightarrow\fal{G_0}\optens\fal{H_0}$.  Thus, via the 
completely isometric injections
\[
\fI\optens\fJ\hookrightarrow\fal{G_0}\optens\fal{H_0}
\cong\fal{G_0\cross H_0}\hookrightarrow\fal{G\cross H}
\]
we may regard $\fI\optens\fJ$ as a compactly supported closed ideal of   
$\fal{G\cross H}$.

Now we have a completely isometric identification
\[
J:\bigl(\ell^1(G)\optens\fI\bigr)\optens\bigl(\ell^1(H)\optens\fJ\bigr)
\to\ell^1(G\cross H)\optens(\fI\optens\fJ).
\]
Using finite sums of elementary tensors we see that  
$q_\fI\otimes q_\fJ=q_{\fI\optens\fJ}\comp J$.
Hence $\ran(q_\fI\otimes q_\fJ)$, with its quotient operator
space structure, must be completely (isometrically) isomorphic
to $\ran q_{\fI\optens\fJ}$, with its quotient operator space
structure.  By projectivity of the operator projective tensor product,
see \cite[7.1.7]{effrosrB}, we have that 
\[
\ran(q_\fI\otimes q_\fJ)=\sofalg\optens\sofalh.
\]
By Corollary \ref{cor:maintheorem1} we have that
$\ran q_{\fI\optens\fJ}=\sofalgh$.   \endpf

If $G$ and $H$ are abelian, the following recovers one of the main
results of Feichtinger \cite{feichtinger}.

\begin{tensorprod1}\label{cor:tensorprod1}
If either $G$ or $H$ admits an open abelian subgroup,
then we have a natural isomorphism $\sofalg\ptens\sofalh
\cong\sofalgh$.
\end{tensorprod1}

\proof If either $G$ or $H$ is abelian, then the proof above can be followed 
almost
verbatim, with the projective tensor product $\ptens$ 
playing the role of the operator projective tensor product $\optens$.  
The reason we require extra hypotheses here is that they
are sufficient (and almost necessary) to obtain that $\falg\ptens\falh\cong
\falgh$, isomorphically, as proved in \cite{losert}, for example.  
Thus we obtain that $\fI\ptens\fJ$ can
be realised as an ideal in $\falgh$.  

If $G$, say, has an open abelian subgroup $A$, with transversal for left cosets
$T$, then by Corollary \ref{cor:maintheorem4} and the reasoning above, 
we obtain isomorphic identifications
\begin{align*}
\sofalg\ptens\sofalh&\cong\ell^1(T)\ptens\sofalb{A}\ptens\sofalh \\
&\cong\ell^1(T\cross\{e_H\})\ptens\sofalb{A\cross H}\cong\sofalgh
\end{align*}
where we obtain the last identification by realising $T\cross\{e_H\}$ 
as a transversal for the left cosets of $A\cross H$ in $G\cross H$.
\endpf

We note that if $G$ and $H$ are both compact, neither having an open abelian
subgroup, then the above result fails by Corollary \ref{cor:maintheorem3} 
(ii) and \cite{losert}.  
We conjecture that our operator space structure on $\sofalg$ is the 
maximal operator space structure
exactly when $G$ admits an open abelian subgroup.  This would imply the result
above.  However, it is clear, only when $G$ has a compact abelian open 
subgroup, that our operator space structure on $\sofalg$ is the maximal one. 
Indeed if $G$ is abelian, then for an arbitrary 
closed ideal $\fI$ of $\falg$, i.e.\ of $\blone{\hat{G}}$, it is not clear
that the subspace operator space structure is the maximal one, whence we
have no means to deduce that $\sofalg\cong\ell^1(G)\optens\fI/\ker q_\fI$
is a maximal operator space.  For results on subspaces of maximal operator 
spaces, see \cite{oikhberg}, for example.

\subsection{Restriction}\label{ssec:restriction}
We recall from \cite{herz,takesakit} that if $H$ is a closed subgroup
of $G$, then the restriction map $u\mapsto u|_H:\falg\to\falh$ is a
quotient map.  In fact, it is a complete quotient map since its adjoint
map is an injective $*$-homomorphism from $\vnh$ onto the von
Neumann algebra generated by $\{\lam(s):s\in H\}$ in $\vng$.

The following result is due to Feichtinger \cite{feichtinger}, in the
abelian case.  However, most of his techniques rely on 
commutativity of $G$, and cannot be adapted to show the general case, even with
no considerations for the operator space structure.

\begin{restrict}\label{theo:restrict}
If $H$ is a closed subgroup in $G$, then the restriction map
\[
u\mapsto u|_H:\sofalg\to\sofalh
\]
is completely surjective.
\end{restrict}

\proof First we must verify that if $u\in\sofalg$, then $u|_H\in\sofalh$.
Let $T$ be a transversal for the right cosets of $H$.  The bijection
$(t,s)\mapsto st:T\cross H\to G$ induces an isomorphism
\[
\ell^1(T)\optens\ell^1(H)=\ell^1(T)\ptens\ell^1(H)
\cong\ell^1(T\cross H)\cong\ell^1(G).
\]  
Now let $\fI=\falkg$, where $K$ is a compact neighbourhood
of the identity in $G$. If $t\in T$, then
\begin{equation}\label{eq:setincl}
(tK)\cap H\subset s_t(K^{-1}K\cap H),\text{ for some }s_t\iin H.
\end{equation}
Indeed, if $(tK)\cap H\not=\varnothing$, then there is $k\in K$
so $tk\in H$, so $t\in Hk^{-1}\subset HK^{-1}$, and thus there is $s_t\iin H$
so $t\in s_tK^{-1}$, whence $tk\in s_tK^{-1}K$.
Now, as in the proof of Lemma \ref{lem:feichdesc}, any $u\iin\sofalg$
can be written in the form
\[
u=\sum_{t\in T}\sum_{s\in H} s\con t\con u_{st},\wwhere
\sum_{t\in T}\sum_{s\in H} \fnorm{u_{st}}<+\infty
\]
and each $u_{st}\in\falkg$.
We then have that for each $t\iin T$, using (\ref{eq:setincl}), that 
\[
\sum_{s\in H} (s\con t\con u_{st})|_H
=\sum_{s\in H} s\con s_t\con\left((s_t^{-1}\con t\con u_{st})|_H\right)
\]
where, $(s_t^{-1}\con t\con u_{st})|_H\in\fK=\falbh{K^{-1}K\cap H}$
and $\fnorm{(s_t^{-1}\con t\con u_{st})|_H}\leq\fnorm{u_{st}}$.
It then follows that $u|_H$, being a 
$\norm{\cdot}_{\ran q_\fK}$-summable series of elements
from $\sofalh$, is itself in $\sofalh$.

Now let us see that restriction is completely surjective.  Let
$\fI$ be as above so that $(\supp\fI)^\circ\cap H\not=\varnothing$.
Note that $\fI|_H$, with the operator space structure given by its
being a quotient of $\fI$ via the restriction map, is a completely
contractive $\falh$-module.  Indeed, this follows from the fact that
$\falh$ is a complete quotient of $\falg$.
Since $\ell^1(H)$ is a (completely) complemented subspace 
of $\ell^1(G)$, we have that $\ell^1(H)\optens\fI$ is a closed
subspace of $\ell^1(G)\optens\fI$.  We have that the following
diagram commutes
\[
\xymatrix{
\tmatinfty{\ell^1(H)\optens\fI}\ar[rrr]^{\tmatinfty{\id\otimes(u\mapsto u|_H)}}
\ar[d]_{\tmatinfty{q_\fI|_{\ell^1(H)\optens\fI}}} 
&  &  & \tmatinfty{\ell^1(H)\optens\fI|_H} \ar[d]^{\tmatinfty{q_{\fI|_H}}} \\
\tmatinfty{\sofalg} \ar[rrr]^{[u_{ij}]\mapsto [u_{ij}|_H]} 
&  &  & \tmatinfty{\sofalh}}
\]
where $q_{\fI|_H}$ is a complete surjection by Corollary
\ref{cor:maintheorem1}, so $\tmatinfty{q_{\fI|_H}}$ is a surjection
by Corollary \ref{cor:tmatop1}.  
Thus the restriction map
$u\mapsto u|_H:\sofalg\to\sofalh$ is completely surjective. \endpf

\subsection{Multipliers on $\wbar{\bltwog}\ptens\bltwog$}
The aim of this section is to develop some techniques for use in the next
section on the averaging operation.  Let 
\[
\tg=\wbar{\bltwog}\ptens\bltwog
\]
where $\wbar{\bltwog}$ denotes the 
conjugate space of $\bltwog$.  We recall that $\tg^*\cong\bdop{\bltwog}$
via the dual pairing $\dpair{\bar{f}\otimes g}{T}=\inprod{Tg}{f}$.
Thus $\tg$ is an operator space with the predual operator space structure.

We may regard $\tg$ as a space of equivalence classes of functions
on $G\cross G$:  if $\ome=\sum_{k=1}^\infty \bar{f}_k\otimes g_k$, where
$\{f_k\}_{k=1}^\infty$ and $\{g_k\}_{k=1}^\infty$ are each summable sequences 
from $\bltwog$, then $\ome(s,t)=\sum_{k=1}^\infty \wbar{f_i(s)}g_i(t)$ for
almost every $s$ and almost every $t$ in $G$.  A function $w:G\cross G\to\Cee$
is called a {\it multiplier} of $\tg$ if for every $\ome\iin\tg$, $m_w\ome$,
defined for almost every $s$ and almost every $t$ in $G$ by 
$m_w\ome(s,t)=w(s,t)\ome(s,t)$, determines an element of $\tg$.  
We let $\mtg$ denote the space of all multipliers $w$ such that 
$m_w:\tg\to\tg$ is a bounded map.  We summarise below, some results
from \cite{spronk}.  We note that there are some trivial differences
between our notations used here, and those in \cite{spronk}; the result
is stated below to be consistent with our present notation.
For a different perspective, we also refer the reader to \cite{neufangrs}.

\begin{tmultipliers}\label{theo:tmultipliers}
{\bf (i)} \cite[Theo.\ 3.3]{spronk}  For each $w\iin\mtg$, $m_w:\tg\to\tg$ is
a completely bounded map with $\norm{m_w}_{\cbop{\tg}}=
\norm{m_w}_{\bdop{\tg}}$.  Thus the space $\mtg\cong\{m_w:w\in\mtg\}$
is a closed subalgebra of $\cbop{\tg}$, and hence a completely contractive
Banach algebra.

{\bf (ii)} \cite[Cor.\ 4.3 \& Theo.\ 5.3]{spronk} If $u\in\fsalg$, the
map $\gam u:G\cross G\to\Cee$, given by $\gam u(s,t)=u(st^{-1})$, 
is an element of $\mtg$.  Moreover, 
$\gam:\fsalg\to\mtg$ is a completely contractive homomorphism.
\end{tmultipliers}

We let $P_G:\tg\to\falg$ be given by
\begin{equation}\label{eq:herzmap}
P_G(\bar{f}\otimes g)=\inprod{\lam(\cdot)g}{f}=\bar{f}\con\check g
\end{equation}
where $\check{g}(s)=g(s^{-1})$ for almost every $s$.  The adjoint map,
$P_G^*:\vng\to\bdop{\bltwog}$, is the inclusion map, hence a complete
isometry, whence $P_G$ is a complete contraction.  We note for $\ome\iin\tg$
that
\[
P_G\ome(s)=\int_G\ome(t,s^{-1}t)dt
\]
for each $s\iin g$.  Thus if $u\in\fsalg$, then
\[
P_G(m_{\gam u}\ome)(s)=\int_Gu\bigl(t(s^{-1}t)^{-1}\bigr)\ome(t,s^{-1}t)dt
=u(s)P_G\ome(s).
\]

We now introduce a class of ideals in $\falg$ which will
prove useful.
Let $K$ be a compact subset of $G$ of positive measure.
Let $\tk=\wbar{\bltwo{K}}\ptens\bltwo{K}$, where we regard
$\bltwo{K}$ as a subspace of $\bltwog$ in the natural way.
We define
\[
\fM(K)=P_G(\tk)
\]
and endow $\fM(K)$ with the quotient operator space structure so it is 
isometrically isomorphic with $\tk/\ker(P_G|_{\tk})$.  Clearly,
$\fM(K)\subset\falg$, and has as a dense subspace $\spn\{\bar{f}\con\check{g}:
f,g\in\bltwo{K}\}$.

\begin{ideals}\label{prop:ideals}
The space $\fM(K)$ is a completely contractive $\fsalg$-module.
Thus it is an ideal in $\falg$ with $\supp\fM(K)\subset K^{-1}K$, and
equipped with an operator space structure by which it is a completely
contractive $\falg$-module.
\end{ideals}

We remark that there is no reason to suspect that $\fM(K)$ is a closed
ideal in $\falg$ for a general compact set $K$.

\medskip
\proof Since $\bltwo{K}$ is a complemented subspace of $\bltwog$,
$\tk$ identifies isometrically as a closed subspace of $\tg$.
As such, $\tk$ is a $\mtg$-submodule of $\tg$.  Also, for $u\iin\fsalg$
and $\ome\iin\tk$ we have $uP_G\ome=P_G(m_{\gam u}\ome$, so
$\fM(K)$ is a $\fsalg$-module.  We thus have that the
following diagram commutes.
\[
\xymatrix{
\fsalg\otimes_\wedge\tk\ar[rr]^{\gam\otimes\id}\ar[d]_{\id\otimes P_G|_{\tk}}
&  & \mtg\optens\tk\ar[rr]^{\quad w\otimes\ome\mapsto m_w\ome}
&  & \tk \ar[d]^{P_G|_{\tk}} \\
\fsalg\otimes_\wedge\fM(K) \ar[rrrr]^{u\otimes v\mapsto uv} &  &  &  & \fM(K)}
\]
Since $\id\otimes P_G|_{\tk}$ is a complete quotient map,
and the maps $\gam\otimes\id$, $w\otimes\ome\mapsto m_w\ome$ 
and $P_G|_{\tk}$ are complete contractions, $u\otimes v\mapsto uv$ 
must be a complete contraction too.

Thus it is clear that $\fM(K)$ is an ideal in $\falg$ and a 
completely contractive $\falg$-module.  It is straightforward to verify
that $\supp\fM(K)\subset K^{-1}K$.  \endpf

Let us note that we may obtain a weak version of a ``tensor product 
factorisation''
result of \cite{feichtinger}.  Let $K$ be a compact subset of $G$
of non-empty interior, and define $q_K^2:\ell^1(G)\ptens\bltwo{K}\to\bltwog$
by $q_K^2(\del_s\otimes f)=s\con f$.  Let $\wtwog=\ran q_K^2$ and norm it
as the quotient space.  It is straightforward to verify that for any other
compact set $K'$, having non-empty interior, that $\ran q_{K'}^2=\ran q_K^2$,
and that the quotient norms are equivalent. 
Let $P_G':\wtwog\ptens\wtwog\to\falg$ be given by
$P_G'(f\otimes g)=f\con\check{g}$.  It can be checked, similarly as in the 
proof of the corollary above, that
$\ran P_G'$, with its quotient norm, is a Segal algebra in $\falg$
on which $G$ acts isometrically by left translations.
Also, $\ran P_G'\subset\sofalg$, and thus, by Theorem  \ref{theo:maintheorem} 
(ii), we obtain $\ran P_G'=\sofalg$.

However, unlike in the commutative case, we do not have that either of
the maps from $\sofalg\ptens\sofalg$ to $\sofalg$, given on elementary tensors
by $u\otimes v\mapsto u\con v$ or $u\otimes v\mapsto u\con\check{v}$,
are surjective.  Indeed, this fails
for compact groups which do not admit an abelian subgroup of finite index,
a fact which follows from \cite[Prop.\ 2.5]{johnson}, in light of Corollary 
\ref{cor:maintheorem3} and the fact that $v\mapsto\check{v}$ is an isometry
on $\falg$.  It would be interesting to know when either of the aforementioned
maps, extended to $\sofalg\optens\sofalg$, surjects onto $\sofalg$.
This fails in general.
Recent work of the author, with B.E. Forrest and E. Samei, has shown that
if $G$ is compact, and hence $\sofalg=\falg$, then each such map
is surjective only when $G$ admits an abelian subgroup of finite index.

\subsection{Averaging over a normal subgroup}\label{ssec:averaging}
Let $N$ be a closed normal subgroup of $G$ and $\tau_N:\bloneg\to\blone{G/N}$
be given for $f\iin\bloneg$ and almost every $sN\in G/N$,
accepting a mild abuse of notation, by
\[
\tau_N(f)(sN)=\int_N f(sn)dn.
\]
This operator is a complete quotient map
as observed in \cite{forrestsw}.  It was shown in \cite{feichtinger},
for an abelian $G$, that $\tau_N(\sofalg)=\sofalb{G/N}$.  We obtain a 
generalisation of that result.

\begin{averaging}\label{theo:averaging}
We have for any locally compact group $G$ with closed normal subgroup $N$
that $\tau_N(\sofalg)=\sofalb{G/N}$, and $\tau_N:\sofalg\to\sofalb{G/N}$
is a complete surjection.
\end{averaging}

\proof We divide the proof into three stages.

{\bf (I)} {\it $\tau_N(\fM(K))\subset\fal{G/N}$ and $\tau_N:\fM(K)\to\fal{G/N}$
is completely bounded.}

Let us first show that $\tau_N(\bltwo{K})\subset\bltwo{G/N}$ and that
$\tau_N:\bltwo{K}\to\bltwo{G/N}$ is bounded.
Let $\vphi:G\to\Cee$ be a continuous function of compact support such that
$\vphi|_K=1$.  Then for any $f\iin\bltwo{K}$ we have, using H\"{o}lder's 
inequality
and the Weyl integral formula, that
\begin{align*}
\bltwonormb{\tau_N(f)}{G/N}^2
&=\int_{G/N}\left|\int_N \vphi(sn)f(sn)dn\right|^2dsN \\
&\leq\int_{G/N}\int_N|\vphi(sn')|^2dn'\int_N|f(sn)|^2dn\, dsN \\
&\leq\sup_{s\in G}\tau_N(|\vphi|^2)(sN)\int_G|f(s)|^2ds.
\end{align*}
We note that $\sup_{s\in G}\tau_N(|\vphi|^2)(sN)<\infty$ since 
$\tau_N(|\vphi|^2)$ is itself continuous and of compact support on $G/N$,
as can be checked using the uniform continuity of $|\vphi|^2$.  
Thus 
\[
\norm{\tau_N}_{\bdop{\bltwo{K},\bltwo{G/N}}}\leq
\sup_{s\in G}\tau_N(|\vphi|^2)(sN)^{1/2}.
\]
Is is shown in \cite[p.\ 187]{lohoue} that for any compactly supported 
$f\in\bloneg$ that
\[
\tau_N(\check{f})=[\Del_{G/N}\tau_N(\check{\Del}_Gf)]^\vee
\]
Now let $\theta_N:\bltwo{K}\to\bltwo{G/N}$ be given by
\[
\theta_N(f)=\Del_{G/N}\tau_N(\check{\Del}_Gf).
\]
Then we have that $\theta_N$ is bounded with
\[
\norm{\theta_N}_{\bdop{\bltwo{K},\bltwo{G/N}}}\leq
\sup_{s\in K}\Del_{G/N}(sN)\norm{\tau_N}_{\bdop{\bltwo{K},\bltwo{G/N}}}
\sup_{t\in K}\frac{1}{\Del_G(t)}.
\]
Now if $f,g\in\bltwo{K}$, then
\[
\tau_N(\bar{f}\con\check{g})=\tau_N(\bar{f})\con\tau_N(\check{g})
=\wbar{\tau_N(f)}\con[\theta_N(g)]^\vee\in\fal{G/N}.
\]
Hence it follows that $\tau_N(\fM(K))\subset\fal{G/N}$.

We now wish to establish that $\tau_N:\fM(K)\to\fal{G/N}$ is completely
bounded.  As in \cite[Sec.\ 3.4]{effrosrB},
we assign $\wbar{\bltwo{K}}$ the row space operator space 
structure, denoted $\wbar{\bltwo{K}}_r$; and $\bltwo{K}$ the column 
operator space structure space, denoted $\bltwo{K}_c$.  Then
we have a completely isometric equality
\[
\tg=\wbar{\bltwo{K}}_r\optens\bltwo{K}_c
\]
by \cite[9.3.2 and 9.3.4]{effrosrB}.  Then, by
\cite[3.4.1 and 7.1.3]{effrosrB},
we have that
\[
\tau_N\otimes\theta_N:\wbar{\bltwo{K}}_r\optens\bltwo{K}_c\to
\wbar{\bltwo{G/N}}_r\optens\bltwo{G/N}_c
\]
is completely bounded, thus is a completely bounded map on $\tk$.
Hence we see that the following diagram commutes
\[
\xymatrix{
\tk\ar[rr]^{\tau_N\otimes\theta_N}
\ar[d]_{P_G|_{\tk}} 
&  & \tb{G/N} \ar[d]^{P_{G/N}} \\
\fM(K) \ar[rr]^{\tau_N} &  & \fal{G/N}}
\]
where $P_G$ is defined in (\ref{eq:herzmap}), and $P_{G/N}$ is defined 
analogously.
Since $P_G|_{\tk}$ is a complete quotient, and $\tau_N\otimes\theta_N$
and $P_{G/N}$ are completely bounded, $\tau_N:\fM(K)\to\fal{G/N}$ 
must be completely bounded too.

{\bf (II)} {\it $\tau_N(\fM(K))$, with its quotient operator space structure,
i.e.\ naturally identified with the quotient space
$\fM(K)/\ker(\tau_N|_{\fM(K)})$, is a completely contractive
$\fal{G/N}$-module.}

Let $\pi_N:G\to G/N$ denote the quotient map.  By \cite[(2.26)]{eymard}
the function $u\mapsto u\comp\pi_N$ defines a complete isometry
from $\fsalgmn$ to $\fsalgn$, the closed subspace of $\fsalg$ of functions
which are constant on cosets of $N$.  Now if $u\in\fsalgmn$ and $v\in\fM(K)$,
then for any $s\iin G$ we have
\[
u(sN)\tau_N(u)(sN)=\int_N u\comp\pi_N(sn)v(sn)dn=\tau_N(u\comp\pi_N\, v)(sN)
\]
Hence $u\tau_N(v)=\tau(u\comp\pi_N\, v)\in\tau_N(\fM(K))$, by Proposition
\ref{prop:ideals}.  

We now establish that $\tau_N(\fM(K))$ is a completely
contractive $\fsalgmn$-module, hence a completely contractive 
$\fal{G/N}$-module.  Letting $\iota:\fsalgn\to\fsalgmn$ be the inverse
of $u\mapsto u\comp\pi_N$, we obtain the following commuting diagram.
\[
\xymatrix{
\fsalgn\otimes_\wedge\fM(K)\ar[rrr]^{u\otimes v\mapsto uv}
\ar[d]_{\iota\otimes\tau_N|_{\fM(K)}} 
&  &  & \fM(K) \ar[d]^{\tau_N|_{\fM(K)}} \\
\fsalgmn\otimes_\wedge\tau_N(\fM(K)) 
\ar[rrr]^{\quad\quad\til{u}\otimes\til{v}\mapsto \til{u}\til{v}} 
&  &  & \fal{G/N}}
\]
Since $\iota\otimes\tau_N|_{\fM(K)}$ is a complete contraction, and 
$u\otimes v\mapsto uv$ and $\tau_N|_{\fM(K)}$ are complete contractions,
we obtain that $\til{u}\otimes\til{v}\mapsto \til{u}\til{v}$
is thus a complete contraction too.

{\bf (III)} {\it The finale.}

We recall from Corollary \ref{cor:maintheorem2} (ii), that the
map $q_{\fM(K)}':\bloneg\optens\fM(K)\to\sofalg$, given on
elementary tensors by $q_{\fM(K)}'f\otimes u=f\con u$
is a complete surjection.  Similarly, appealing also to (II) above,
we have that $q_{\tau_N(\fM(K))}':\blone{G/N}\optens\tau_N(\fM(K))\to
\sofalb{G/N}$ is a complete surjection.
It is then clear that the following diagram commutes.
\[
\xymatrix{
\bloneg\optens\fM(K)\ar[rr]^{\tau_N\otimes\tau_N}
\ar[d]_{q_{\fM(K)}'} 
&  & \blone{G/N}\optens\tau_N(\fM(K)) \ar[d]^{q_{\tau_N(\fM(K))}'} \\
\sofalg \ar[rr]^{\tau_N} &  & \blone{G/N}}
\]
In particular $\tau_N(\sofalg)\subset\ran q_{\tau_N(\fM(K))}'=\sofalb{G/N}$. 
Moreover, since $q_{\fM(K)}'$ and $\tau_N\otimes\tau_N$, as above, and
$q_{\tau_N(\fM(K))}':\blone{G/N}\optens\tau_N(\fM(K))\to\sofalb{G/N}$ 
are complete surjections,
$\tau_N:\sofalg\to\sofalb{G/N}$ is completely bounded.  
Amplifying the diagram by $\stmatinfty$, and appealing
to Corollary \ref{cor:tmatop1}, as in the end of the proof of Theorem
\ref{theo:restrict}, we see that
$\tau_N:\sofalg\to\sofalb{G/N}$ is completely surjective.  \endpf

We remark that it follows from the above theorem that $\tau_N(\falcg)
\subset\mathrm{A}_c(G/N)$.  We note that it was proved in \cite{forrestsw}
that $\tau_N(\lebfalg)=\blone{G/N}$.  We also note that it was shown by
Lohou\'{e} \cite{lohoue} that $\tau_N:\falkg\to\fal{G/N}$ is a bounded map, with
bounded depending on $K$.  This fact can be deduced from our result, but
Lohou\'{e}'s proof is much simpler, though it is not obvious how to adapt
his proof to show that $\tau_N\in\cbop{\falkg,\fal{G/N}}$.


\subsection{An Isomorphism Theorem}\label{ssec:isomorphism}
Let $G$ and $H$ be locally compact groups.  Wendel \cite{wendel} proved
that there is an isometric isomorphism between the convolution algebras
$\bloneg$ and $\blone{H}$ if and only if $G$ and $H$ are isomorphic
topological groups \cite{wendel}.  Also, Walter \cite{walter} proved that
$\falg$ and $\falh$ are isometrically isomorphic if and only if $G$ and $H$ 
are isomorphic topological groups.  Since we lack fixed norms
on our algebras $\sofalg$ and $\sofalh$, it is not reasonable to expect
and ``isometric isomorphism'' theorem, in the spirit of Wendel's and Walter's
theorems.  In fact, if $G$ and $H$ are both discrete groups
having the same cardinality, then Corollary \ref{cor:maintheorem3} tells
us there is a multiplicative isomorphism identifying $\sofalg\cong\sofalh$.
Similarly, if $G$ and $H$ are finite abelian groups, then there is a 
convolutive isomorphism identifying $\sofalg\cong\sofalh$.
To obtain a satisfactory result, we must simultaneously
exploit the facts $\sofalg$ and $\sofalh$ are pointwise and convolutive 
algebras.

\begin{isomorphism}\label{theo:isomorphism}
Let $G$ and $H$ be locally compact groups and $\Phi:\sofalg\to\sofalh$
be a bounded linear bijection which satisfies
\[
\Phi(uv)=\Phi u\,\Phi v \quad\aand\quad
\Phi(u\con v)=\Phi u\con\Phi v
\]
for every $u,v\iin\sofalg$.  There is a homeomorphic isomorphism
$\alp:G\to H$ such that
\[
\Phi u=u\comp\alp
\]
for each $u\iin\sofalg$.
\end{isomorphism}

\proof We recall that $\falg$ has Gelfand spectrum $G$, implemented by
evaluation functionals.  Since $\sofalg$ is a Segal algebra in $\falg$,
it follows from \cite[Thm.\ 2.1]{burnham} that $\sofalg$ has Gelfand spectrum
$G$ too.  The same holds for $\sofalh$.  Thus we may define $\alp:H\to G$
by letting for $h\iin H$, $\alp(h)$ be the element of $G$ which satisfies
$u\bigl(\alp(h)\bigr)=\Phi u(h)$ for each $u\in\sofalg$.  Then $\alp$ is 
continuous.  Indeed, if not, we may find a net  $h_i\to h$ in $H$ and
a neighbourhood $U$ of $\alp(h)$, such that
$\alp(h_i)\not\in U$ for each $i$.  Using regularity we may find
$u\iin\sofalg$ such that $u\bigl(\alp(h)\bigr)=1$ and $\supp u\subset U$.
But then we would obtain 
\[
\lim_i \Phi u(h_i)=\lim_i u\bigl(\alp(h_i)\bigr)=0\not=1=u\bigl(\alp(h)\bigr)
=\Phi u(h)
\]
which contradicts that $\Phi u$ is continuous, in particular that
$\Phi u\in\sofalh$.  We may similarly obtain a continuous map $\beta:G\to H$
satisfying $v\bigl(\beta(s)\bigr)=\Phi^{-1}v(s)$ for all $s\in G$ and
$v\iin\sofalh$. We clearly have that $\beta=\alp^{-1}$, hence $\alp$ is a 
homeomorphism.

It remains to see that $\alp$ is a group homomorphism.  Let $\fU$ be 
a compact neighbourhood basis of the identity $e_G$ in $G$.  For each
$U\iin\fU$ find $u_U\iin\sofalg$ such that
\[
\supp u_U\subset U\quad\aand\quad
\int_G |u_U(s)|ds=1.
\]
Then $(u_U)$ is a bounded approximate identity for $\bloneg$,
hence a convolutive approximate unit for $\sofalg$.  Since $\Phi$ is
a surjective convolutive homomorphism, $(\Phi u_U)$ is a convolutive
approximate identity for $\sofalh$.  Let $h_1,h_2\iin H$ and suppose
that $\alp(h_1)^{-1}\alp(h_2)\not=\alp(h_1^{-1}h_2)$.  We could then
find $v\in\sofalg$ such that
\[
v\bigl(\alp(h_1)^{-1}\alp(h_2)\bigr)=1\aand
v(s)=0\text{ for all }s\text{ in a nbhd.\ of } \alp(h_1^{-1}h_2).
\]
Then we would have that
\begin{align*}
v\bigl(\alp(h_1)^{-1}\alp(h_2)\bigr)
&=\alp(h_1)\con v\bigl(\alp(h_2)\bigr) 
=\lim_{U}\alp(h_1)\con u_U\con v\bigl(\alp(h_2)\bigr) \\
&=\lim_{U}\Phi\bigl(\alp(h_1)\con u_U\con v\bigr)(h_2) \\
&=\lim_{U}\Phi\bigl(\alp(h_1)\con u_U\bigr)\con\Phi v(h_2) \\
&=\lim_{U}\int_H u_U\bigl(\alp(h_1)^{-1}\alp(r)\bigr)
v\bigl(\alp(r^{-1}h_2)\bigr)dr \\
&=0
\end{align*}
where we obtain the last equality from the fact that
\[
\supp\bigl(\alp(h_1)\con u\bigr)\comp\alp
\subset\{r\in H:\alp(h_1)^{-1}\alp(r)\in U\}
=\alp^{-1}\bigl(\alp(h_1)U\bigr)
\]
and the supposition that $v\bigl(\alp(r^{-1}h_2)\bigr)=0$ for all 
$r\iin\alp^{-1}\bigl(\alp(h_1)U\bigr)$, for a suitably small choice of $U$.
This contradicts that $v\bigl(\alp(h_1)^{-1}\alp(h_2)\bigr)=1$, whence such
a $v$ cannot be chosen, and we thus conclude that $\alp(h_1)^{-1}\alp(h_2)
=\alp(h_1^{-1}h_2)$.  Substituting $e_H=h_1$, we see that $\alp(e_H)=e_G$ and
then, substituting $e_H$ for $h_2$, we obtain 
that $\alp(h_1^{-1})=\alp(h_1)^{-1}$ for each $h_1\iin H$.  Thus
$\alp$ is a group homomorphism.  \endpf

Our theorem above is not special to the class of algebras $\sofalg$.  
In fact it
can be applied to any class of regular Banach algebras with spectrum  $G$,
each of which is a Segal algebra of $\bloneg$.
Examples of such are $\lebfalg$ from (\ref{eq:lebfalg}) and the Wiener
algebra $\mathrm{W}_0(G)$ as defined in \cite{feichtinger81}.

{
\bibliography{segalfourierbib}

\begin{thebibliography}{10}

\bibitem{blecher}
D.~P. Blecher.
\newblock The standard dual of an operator space.
\newblock {\em Pacific Math. J.}, 153(1):15--30, 1992.

\bibitem{burnham}
J.~T. Burnham.
\newblock Closed ideals in subalgebras of {B}anach algebras, {I}.
\newblock {\em Proc. Amer. Math. Soc.}, 32:551--555, 1972.

\bibitem{effrosr}
E.~G. Effros and Z.-J. Ruan.
\newblock On approximation properties for operator spaces.
\newblock {\em International J. Math.}, 1:163--187, 1990.

\bibitem{effrosrB}
E.~G. Effros and Z.-J. Ruan.
\newblock {\em Operator Spaces}, volume~23 of {\em London Math. Soc., New
  Series}.
\newblock Claredon Press, Oxford Univ. Press, New York, 2000.

\bibitem{eymard}
P.~Eymard.
\newblock L'algebre de {F}ourier d'un groupe localement compact.
\newblock {\em Bull. Soc. Math. France}, 92:181--236, 1964.

\bibitem{feichtinger81}
H.~G. Feichtinger.
\newblock A characterization of minimal homogeneous {B}anach spaces.
\newblock {\em Proc. Amer. Math. Soc.}, 81:55--61, 1981.

\bibitem{feichtinger}
H.~G. Feichtinger.
\newblock On a new {S}egal algebra.
\newblock {\em Montash. Math.}, 92:269--289, 1981.

\bibitem{feichtingerg}
H.~G. Feichtinger and K.~Gr\"{o}chening.
\newblock Gabor frames and time-frequency analysis of distributions.
\newblock {\em J. Funct. Anal.}, 146:464--495, 1997.

\bibitem{forrestsw}
B.~E. Forrest, N.~Spronk, and P.~J. Wood.
\newblock Operator {S}egal algebras in {F}ourier algebras.
\newblock Accepted in {\it Studia Math.}, see ArXiv {\tt math.FA/0601509},
  2007.

\bibitem{ghahramanil}
F.~Ghahramani and A.~T.~M. Lau.
\newblock Weak amenability of certain classes of {B}anach algebras without
  bounded approximate units.
\newblock {\em Math. Proc. Camb. Phil. Soc.}, 133:357--371, 2002.

\bibitem{ghahramanil2}
F.~Ghahramani and A.~T.~M. Lau.
\newblock Approximate weak amenability, derivations and {A}rens regularity of
  {S}egal algebras.
\newblock {\em Stud. Math}, 169:189--205, 2005.

\bibitem{herz}
C.~S. Herz.
\newblock Harmonic synthesis for subgroups.
\newblock {\em Ann. Inst. Fourier (Grenoble)}, 23(3):91--123, 1973.

\bibitem{hewittrII}
E.~Hewitt and K.~A. Ross.
\newblock {\em Abstract Harmonic Analysis II}, volume 152 of {\em Die
  Grundlehren der mathematischen Wissenschaften}.
\newblock Springer, Berlin Hieldelberg, 1970.

\bibitem{hewittrI}
E.~Hewitt and K.~A. Ross.
\newblock {\em Abstract Harmonic Analysis I}, volume 115 of {\em Grundlehern
  der mathemarischen Wissenschaften}.
\newblock Springer, New York, second edition, 1979.

\bibitem{johnsonm}
B.~E. Johnson.
\newblock {\em Cohomology in Banach Algebras}, volume 127 of {\em Mem. Amer.
  Math. Soc.}
\newblock Amer. Math. Soc., Providence, RI, 1972.

\bibitem{johnson}
B.~E. Johnson.
\newblock Non-amenability of the {F}ourier algebra of a compact group.
\newblock {\em J. London Math. Soc.}, 50(2):361--374, 1994.

\bibitem{lohoue}
N.~Lohou\'{e}.
\newblock Remarques sur les ensembles de synth\`{e}se des alg\`{e}bras de
  groupe localement compact.
\newblock {\em J. Funct. Anal.}, 13:185--194, 1973.

\bibitem{losert}
V.~Losert.
\newblock On tensor products of the {F}ourier algebras.
\newblock {\em Arch. Math.}, 43:370--372, 1984.

\bibitem{neufangrs}
M.~Neufang, Z.-J. Ruan, and N.~Spronk.
\newblock Completely isometric representations of ${M}_{cb}{A}({G})$ and
  ${UCB}(\hat{{G}})$.
\newblock Accepted in {\it Trans. Amer Math. Soc.}, 2006.

\bibitem{oikhberg}
T.~Oikhberg.
\newblock Subspaces of maximal operator spaces.
\newblock {\em Integr. Equ. Oper. Theory}, 48:81--102, 2004.

\bibitem{paterson}
A.L.T. Paterson.
\newblock {The class of locally compact groups $G$ for which $C^*(G)$ is
  amenable}.
\newblock In {\em Harmonic analysis (Luxembourg, 1987)}, pages 226--237,
  Berlin, 1988. Springer.

\bibitem{reiter}
H.~Reiter.
\newblock {\em Classical Harmonic Analysis and Locally Compact Groups}.
\newblock Oxford Math. Monographs. Claredon, Oxford, 1968.

\bibitem{reiters}
H.~Reiter and J.~D. Stegeman.
\newblock {\em Classical Harmonic Analysis and Locally Compact Groups},
  volume~22 of {\em London Math. Soc., New Series}.
\newblock Claredon Press, Oxford Univ. Press, New York, 2000.

\bibitem{spronk}
N.~Spronk.
\newblock Measurable {S}chur multipliers and completely bounded multipliers of
  the {F}ourier algebras.
\newblock {\em Proc. London Math. Soc.}, 89:161--192, 2004.

\bibitem{takesakit}
M.~Takesaki and N.~Tatsuuma.
\newblock Duality and subgroups, {II}.
\newblock {\em J. Funct. Anal.}, 11:184--190, 1972.

\bibitem{walter}
M.~E. Walter.
\newblock ${W}\sp{*} $-algebras and nonabelian harmonic analysis.
\newblock {\em J. Funct. Anal.}, 11:17--38, 1972.

\bibitem{wendel}
J.~G. Wendel.
\newblock Left centralizers and isomorphisms of group algebras.
\newblock {\em Pacific J. Math.}, 2:251--261, 1952.

\end{thebibliography}
\bibliographystyle{plain}
}

\end{document}